\documentclass[10pt]{article}
\usepackage{amssymb,t1enc}
\newcommand{\wt}{\widetilde}
\newcommand{\R}{\mathbb R}
\newcommand{\C}{\mathbb C}

\newcommand{\F}{\mathbb F}
\newcommand{\A}{\textbf{\textit{A}}}
\newcommand{\B}{\textbf{\textit{B}}}
\newcommand{\ii}{\textbf{\textit{i}}}
\begin{document}
\footnotetext{
\footnotesize
\par
\noindent
{\it 2000 Mathematics Subject Classification:} 39B32, 51M05.
\\
{\it Key words and phrases:} affine (semi-affine) isometry,
affine (semi-affine) mapping with orthogonal linear part,
Beckman-Quarles theorem,
unit-distance preserving mapping.}
\par
\noindent
\centerline{\Large The Beckman-Quarles theorem for continuous}
\par
\noindent
\centerline{\Large mappings from ${\C}^n$ to ${\C}^n$}
\vskip 0.8truecm
\par
\centerline{\Large Apoloniusz Tyszka}
\vskip 0.8truecm
\par
\noindent
{\bf Summary.} Let $\varphi_n:{\C}^n \times {\C}^n \to \C$,
$\varphi_n((x_1,...,x_n),(y_1,...,y_n))=\sum_{i=1}^n\limits (x_i-y_i)^2$.
We say that $f:{\C}^n \to {\C}^n$ preserves distance $d \in \C$,
if for each $X,Y \in~{\C}^n$ $\varphi_n(X,Y)=d^2$ implies
$\varphi_n(f(X),f(Y))=d^2$. We prove:
if $n \geq 2$ and a continuous $f:{\C}^n \to {\C}^n$
preserves unit distance, then $f$ has a form
$I \circ (\underbrace{\rho,...,\rho}_{n-{\rm times}})$,
where $I:{\C}^n \to {\C}^n$ is an affine
mapping with orthogonal linear part and $\rho:~\C \to \C$
is the identity or the complex conjugation.
For $n \geq 3$ and bijective~$f$
the theorem follows from Theorem~2 in \cite{Schroder}.
\vskip 0.8truecm
\par
The classical Beckman-Quarles theorem states that each
unit-distance preserving mapping from
${\R}^n$ to ${\R}^n$ ($n \geq 2$) is an isometry,
see \cite{Beckman-Quarles}-\cite{Everling}.
Let $\varphi_n:{\C}^n \times {\C}^n \to \C$,
$\varphi_n((x_1,...,x_n),(y_1,...,y_n))=\sum_{i=1}^n\limits (x_i-y_i)^2$.
We say that $f:{\C}^n \to {\C}^n$ preserves distance $d \in \C$,
if for each $X,Y \in {\C}^n$ $\varphi_n(X,Y)=d^2$ implies
$\varphi_n(f(X),f(Y))=d^2$.
If $f:{\C}^n \to {\C}^n$ and for each $X,Y \in {\C}^n$
$\varphi_n(X,Y)=\varphi_n(f(X),f(Y))$, then $f$ is an affine mapping
with orthogonal linear part;
it follows from a general theorem proved in [3,~58~ff], see also [4, p. 30].
\vskip 0.2truecm
\par
\noindent
Let $D({\R}^n,{\C}^n)$ denote the set
of all positive numbers $d$ with the property:
\vskip 0.2truecm
\par
if $X,Y \in {\R}^n$ and $\varphi_n(X,Y)=d^2$,
then there exists a finite set $S_{XY}$ with
$\left \{X,Y \right\} \subseteq S_{XY} \subseteq {\R}^n$
such that any map $f:S_{XY}\rightarrow {\C}^n$ that preserves
unit distance satisfies $\varphi_n(X,Y)=\varphi_n(f(X),f(Y))$.
\vskip 0.2truecm
\par
\noindent
Obviously, if $d \in D({\R}^n,{\C}^n)$ and $f:{\R}^n \to {\C}^n$
preserves unit distance, then $f$ preserves distance $d$.
\vskip 0.2truecm
\par
\noindent
Let $D({\C}^n,{\C}^n)$ denote the set of all positive numbers $d$
with the property:
\vskip 0.2truecm
\par
if $X,Y \in {\C}^n$ and $\varphi_n(X,Y)=d^2$,
then there exists a finite set $S_{XY}$ with
$\left\{X,Y \right\} \subseteq S_{XY} \subseteq {\C}^n$ such
that any map $f:S_{XY}\rightarrow {\C}^n$ that preserves
unit distance satisfies $\varphi_n(X,Y)=\varphi_n(f(X),f(Y))$.
\vskip 0.2truecm
\par
\noindent
Obviously, if $d \in D({\C}^n,{\C}^n)$ and $f:{\C}^n \to {\C}^n$
preserves unit distance, then $f$ preserves distance $d$.
\vskip 0.2truecm
\par
If $d>0$, $X,Y \in {\C}^n$ and $\varphi_n(X,Y)=d^2$,
then there exists an affine $I:{\C}^n \to {\C}^n$ such that
$I((0,\underbrace{0,~...,~0}_{n-1~{\rm times}}))=X$, $I((d,\underbrace{0,~...,~0}_{n-1~{\rm times}}))=Y$,
and a linear part of $I$ is orthogonal. Hence,
\vskip 0.2truecm
\par
\noindent
{\bf (1)}~~$D({\R}^n,{\C}^n) \subseteq D({\C}^n,{\C}^n)$.
\vskip 0.2truecm
\par
\noindent
The author proved in \cite{Tyszka2}:
\vskip 0.2truecm
\par
\noindent
{\bf (2)}~~$D({\R}^n,{\C}^n)$ is a dense subset of $(0,\infty)$ for each $n \geq 2$.
\vskip 0.2truecm
\par
\noindent
From {\bf (2)} we obtain (\cite{Tyszka2}):
\vskip 0.2truecm
\par
\noindent
{\bf (3)}~~if $n \geq 2$ and a continuous $f:{\R}^n \to {\C}^n$ preserves unit
distance, then $f$ preserves all positive distances.
\vskip 0.2truecm
\par
\noindent
Obviously,
\vskip 0.2truecm
\par
\noindent
{\bf (4)}~~if $n \geq 1$ and $f:{\R}^n \to {\C}^n$
preserves all positive distances, then there exists an
affine $I:{\C}^n \to {\C}^n$ such that a linear
part of $I$ is orthogonal and $I_{|{\R}^{\scriptstyle n}}=f$.
\vskip 0.2truecm
\par
\noindent
By {\bf (1)} and {\bf (2)}:
\vskip 0.2truecm
\par
\noindent
{\bf (5)}~~$D({\C}^n,{\C}^n)$ is a dense subset of $(0,\infty)$ for each $n \geq 2$.
\vskip 0.2truecm
\par
\noindent
As a corollary of {\bf (5)} we obtain (cf. {\bf (3)}):
\vskip 0.2truecm
\par
\noindent
{\bf Lemma 1.} If $n \geq 2$ and a continuous $f:{\C}^n \to {\C}^n$
preserves unit distance, then $f$ preserves all positive distances.
\vskip 0.2truecm
\par
\noindent
{\it Proof.} Let $d>0$, $X,Y \in {\C}^n$, $\varphi_n(X,Y)=d^2$.
By {\bf (5)} there exists a sequence $\{d_m\}_{m=1,2,3,...}$ tending
to $d$, where all $d_m$ belong to $D({\C}^n,{\C}^n)$. Since
$f$ and $\varphi_n$ are continuous, and $f$ preserves all distances
$d_m$ ($m=1,2,3,...$),
\vskip 0.2truecm
\par
\noindent
\begin{eqnarray*}
\varphi_n(f(X),f(Y))=
\lim_{m\rightarrow\infty}\limits \varphi_n\left(f\left(\frac{d_m}{d}X\right),f\left(\frac{d_m}{d}Y\right)\right)&=&\\
\lim_{m\rightarrow\infty}\limits \varphi_n\left(\frac{d_m}{d}X,\frac{d_m}{d}Y\right)=
\lim_{m\rightarrow\infty}\limits d_m^2=d^2.
\end{eqnarray*}
\vskip 0.2truecm
\par
\noindent
Let $\tau:\C \to \C$ denote the complex conjugation.
\vskip 0.2truecm
\par
\noindent
{\bf Theorem 1.} If $n \geq 2$, $f:{\C}^n \to {\C}^n$
preserves unit distance and $f_{|{\R}^{\scriptstyle n}}={\rm id}({\R}^n)$,
then $f(X) \in \{X, (\underbrace{\tau,...,\tau}_{n-{\rm times}})(X)\}$
for each $X \in {\C}^n$.
\vskip 0.2truecm
\par
\noindent
{\it Proof.} It is true if $X \in {\R}^n$.
Assume now that $X=(x_1,...,x_n) \in {\C}^n \setminus {\R}^n$.
Let $x_1=a_1+b_1 \cdot \ii$,~...,~$x_n=a_n+b_n \cdot \ii$, where
$a_1,b_1$, ..., $a_n,b_n \in \R$. We choose $j \in \{1,...,n\}$
with $b_j \neq 0$. For $k \in \{1,...,n\} \setminus \{j\}$ and $t \in \R$
we define $S_k(t)=(s_{k,1}(t),...,s_{k,n}(t)) \in {\R}^n$ as follows:
\par
\noindent
$s_{k,j}(t)=a_j+tb_k$,
\par
\noindent
$s_{k,k}(t)=a_k-tb_j$,
\par
\noindent
$s_{k,i}(t)=a_i$, if $i \in \{1,...,n\} \setminus \{j,k\}$.
\vskip 0.2truecm
\par
\noindent
Let $t_0=\sqrt{\frac{\textstyle \sum_{i=1}^n\limits b_i^2}{\textstyle b_j^2}}$.
For each $k \in \{1,...,n\} \setminus \{j\}$ and all $t \in \R$
\par
\noindent
\centerline{$\varphi_n((x_1,...,x_n),S_k(t))=t^2(b_j^2+b_k^2) - \sum_{i=1}^n\limits b_i^2$.}
\par
\noindent
By this, for each $k \in \{1,...,n\} \setminus \{j\}$ and all $t \geq t_0$
$$\varphi_n((x_1,...,x_n),S_k(t)) \geq 0.$$
By {\bf (5)}:
\vskip 0.2truecm
\par
\noindent
{\bf (6)}~~for each $k \in \{1,...,n\} \setminus \{j\}$ the set
$$A_k(x_1,...,x_n):=\{t \in \R: \varphi_n((x_1,...,x_n),S_k(t)) \in D({\C}^n,{\C}^n)\}$$
is a dense subset of $(t_0,\infty)$.
\vskip 0.2truecm
\par
\noindent
Let $f((x_1,...,x_n))=(y_1,...,y_n)$. Since $f$ preserves all distances
in $D({\C}^n,{\C}^n)$, for each
$k \in \{1,...,n\} \setminus \{j\}$ and all $t \in A_k(x_1,...,x_n)$
$$
\varphi_n((x_1,...,x_n),S_k(t))=
\varphi_n(f((x_1,...,x_n)),f(S_k(t)))=
\varphi_n((y_1,...,y_n),S_k(t)).
$$
\par
\noindent
Hence, for each $k \in \{1,...,n\}\setminus \{j\}$ and all
$t \in A_k(x_1,...,x_n)$
$$
t^2(b_j^2+b_k^2)-\sum_{i=1}^n b_i^2=
(y_{j}-a_{j}-tb_{k})^2+(y_k-a_k+tb_j)^2+
\sum_{i \in \{1,...,n\} \setminus \{j,k\}} (y_i-a_i)^2.
$$
\vskip 0.2truecm
\par
\noindent
Thus, for each $k \in \{1,...,n\} \setminus \{j\}$ and all
$t \in A_k(x_1,...,x_n)$
$$
\sum_{i=1}^n (y_i-a_i)^2+\sum_{i=1}^n b_i^2=2t \cdot (b_k(y_j-a_j)-b_j(y_k-a_k)).
$$
\vskip 0.2truecm
\par
\noindent
Hence by {\bf (6)}:
\vskip 0.2truecm
\par
\noindent
{\bf (7)}~~$y_k-a_k=\frac{\textstyle b_k}{\textstyle b_j} \cdot (y_j-a_j)$
for each $k \in \{1,...,n\} \setminus \{j\}$
\vskip 0.2truecm
\par
\noindent
and
\vskip 0.2truecm
\par
\noindent
{\bf (8)}~~$\sum_{i=1}^n\limits (y_i-a_i)^2+\sum_{i=1}^n\limits b_i^2=0$.
\vskip 0.2truecm
\par
\noindent
Applying {\bf (7)} to {\bf (8)} we obtain
$$
(y_j-a_j)^2 + \sum_{k \in \{1,...,n\} \setminus \{j\}} \frac{b_k^2}{b_j^2} \cdot (y_j-a_j)^2 + \sum_{i=1}^n b_i^2=0.
$$
It gives $\left(\frac{\textstyle (y_j-a_j)^2}{\textstyle b_j^2}+1 \right) \cdot \sum_{i=1}^n\limits b_i^2=0$.
Since $\sum_{i=1}^n\limits b_i^2 \neq 0$, we get
$$\underbrace{y_j=a_j+b_j \cdot \ii=x_j}_{\rm case~1}~~{\rm or}~~\underbrace{y_j=a_j-b_j \cdot \ii=\tau(x_j)}_{\rm case~2}.$$
\par
\noindent
In case 1, by {\bf (7)} for each $k\in \{1,...,n\} \setminus \{j\}$
\vskip 0.2truecm
\par
\noindent
\centerline{$y_k=a_k+\frac{\textstyle b_k}{\textstyle b_j} \cdot (y_j-a_j)=a_k+\frac{\textstyle b_k}{\textstyle b_j} \cdot (a_j+b_j \cdot \ii -a_j)=a_k+b_k \cdot \ii=x_k$.}
\vskip 0.2truecm
\par
\noindent
In case 2, by {\bf (7)} for each $k \in \{1,...,n\} \setminus \{j\}$
\vskip 0.2truecm
\par
\noindent
\centerline{$y_k=a_k+\frac{\textstyle b_k}{\textstyle b_j} \cdot (y_j-a_j)=a_k+\frac{\textstyle b_k}{\textstyle b_j} \cdot (a_j-b_j \cdot \ii -a_j)=a_k-b_k \cdot \ii=\tau(x_k)$.}
\vskip 0.2truecm
\par
\noindent
The proof is completed.
\vskip 0.2truecm
\par
Let $n \geq 2$, $f:{\C}^n \to {\C}^n$ preserves unit distance,
$f_{|{\R}^{\scriptstyle n}}={\rm id}({\R}^n)$. We would like to prove:
$f={\rm id}({\C}^n)$ or $f=(\underbrace{\tau,...,\tau}_{n-{\rm times}})$;
we will prove it later in Theorem 2.
By Theorem 1 the sets
\vskip 0.2truecm
\par
\noindent
~~~~~~~~~~~~~~~~~~~~~~~~~~~~$\A=\{X \in {\C}^n: f(X)=X\}$
\par
\noindent
and
\par
\noindent
~~~~~~~~~~~~~~~~~~~~~~~~~~~~$\B=\{X \in {\C}^n: f(X)=(\underbrace{\tau,...,\tau}_{n-{\rm times}})(X)\}$
\vskip 0.2truecm
\par
\noindent
satisfy $\A \cup \B={\C}^n$.
\vskip 0.2truecm
\par
\noindent
Let $\psi_n: {\C}^n \times {\C}^n \to \R$, $\psi_n((x_1,...,x_n),(y_1,...,y_n))=\sum_{k=1}^n\limits{\rm Im}(x_k) \cdot {\rm Im}(y_k)$.
\vskip 0.2truecm
\par
\noindent
{\bf Lemma 2.} If $x_1,...,x_n,y_1,...,y_n \in \C$,
$\varphi_n((x_1,...,x_n),(y_1,...,y_n))=1$ and
\par
\noindent
$\psi_n((x_1,...,x_n),(y_1,...,y_n)) \neq 0$, then
\vskip 0.2truecm
\par
\noindent
{\bf (9)}~~~~~~~~~~~~~~~~~~~~~~~~~$(y_1,...,y_n) \in \A$ implies $(x_1,...,x_n) \in \A$
\par
\noindent
and
\par
\noindent
{\bf (10)}~~~~~~~~~~~~~~~~~~~~~~~~$(y_1,...,y_n) \in \B$ implies $(x_1,...,x_n) \in \B$.
\vskip 0.2truecm
\par
\noindent
{\it Proof.} We prove only {\bf (9)}, the proof of {\bf (10)}
follows analogically. Assume, on the contrary, that $(y_1,...,y_n) \in \A$
and $(x_1,...,x_n) \not\in \A$. Since $\A \cup \B={\C}^n$,
$(x_1,...,x_n) \in \B$.
Let $x_1=a_1+{b_1} \cdot \ii$,~...,~$x_n=a_n+{b_n} \cdot \ii$,
$y_1=\wt{a}_1+\wt{b}_1 \cdot \ii$,~...,~$y_n=\wt{a}_n+\wt{b}_n \cdot \ii$,
where $a_1,b_1,...,a_n,b_n,\wt{a}_1,\wt{b}_1,...,\wt{a}_n,\wt{b}_n \in \R$.
Since $f$ preserves unit distance,
\vskip 0.2truecm
\par
\noindent
{\bf (11)}~~~~~~~~~~~~~~~$\varphi_n((x_1,...,x_n),(y_1,...,y_n))=$
\vskip 0.2truecm
~~~~~~~~$\varphi_n(f((x_1,...,x_n)),f((y_1,...,y_n)))=
\sum_{k=1}^n\limits (a_k-b_k \cdot \ii-\wt{a}_k-\wt{b}_k \cdot \ii)^2.$
\vskip 0.2truecm
\par
\noindent
Obviously,
\par
\noindent
{\bf (12)}~~~~~~~~~~~~~~~$\varphi_n((x_1,...,x_n),(y_1,...,y_n))=\sum_{k=1}^n\limits (a_k+b_k \cdot \ii - \wt{a}_k-\wt{b}_k \cdot \ii)^2$.
\vskip 0.2truecm
\par
\noindent
Subtracting {\bf (11)} and {\bf (12)} by sides we obtain
$$4 \sum_{k=1}^n b_k \wt{b}_k+ 4\sum_{k=1}^n b_k(a_k-\wt{a}_k) \cdot \ii=0,$$
so in particular $\psi_n((x_1,...,x_n),(y_1,...,y_n))=\sum_{k=1}^n\limits b_k\wt{b}_k=0,$
a contradiction.
\vskip 0.2truecm
\par
\noindent
The next lemma is obvious.
\vskip 0.2truecm
\par
\noindent
{\bf Lemma 3.} For each $S,T \in {\R}^n$ there exist $m \in \{1,2,3,...\}$
and $P_1,...,P_m \in {\R}^n$ such that
$||S-P_1||=||P_1-P_2||=...=||P_{m-1}-P_m||=||P_m-T||=1$.
\vskip 0.2truecm
\par
\noindent
{\bf Lemma 4.} For each $X \in {\C}^n \setminus {\R}^n$
\vskip 0.2truecm
\par
\noindent
\centerline{$(\underbrace{\ii,...,\ii}_{n-{\rm times}}) \in \A$ implies $X \in \A$}
\par
\noindent
and
\par
\noindent
\centerline{$(\underbrace{\ii,...,\ii}_{n-{\rm times}}) \in \B$ implies $X \in \B$.}
\vskip 0.2truecm
\par
\noindent
{\it Proof.} Let $X=(a_1+b_1 \cdot \ii,~...,~a_n+b_n \cdot \ii)$, where
$a_1,...,a_n,b_1,...,b_n \in \R$. We choose $j \in \{1,...,n\}$ with
$b_j \neq 0$. The points
\vskip 0.2truecm
\par
\noindent
\begin{eqnarray*}
S=\Biggl(a_{1}+\sqrt{\frac{1+(b_{j}-1)^2}{n-1}},~...,~
a_{j-1}+\sqrt{\frac{1+(b_{j}-1)^2}{n-1}}~,& &\\
\underbrace
{a_{j}+\sqrt{1+\sum_{{k \in \{1,...,n\}\setminus\{j\}}}\limits b_k^2}}
_{j-{\rm th}~{\rm coordinate}}~,& &\\
a_{j+1}+\sqrt{\frac{1+(b_{j}-1)^{2}}{n-1}},~...,~
a_{n}+\sqrt{\frac{1+(b_{j}-1)^2}{n-1}}\Biggr)
\end{eqnarray*}
\vskip 0.2truecm
\par
\noindent
and $T=(0,~...,~0,\underbrace{~~~~~~\sqrt{n}~~~~~~}_{j-{\rm th}~{\rm coordinate}},0,~...,~0)$ belong to ${\R}^n$.
Applying Lemma 3 we find
$m \in \{1,2,3,...\}$
and
$P_1,...,P_m \in {\R}^n$ satisfying $||S-P_1||=||P_1-P_2||= ... =||P_{m-1}-P_m||=||P_m-T||=1$.
The points
\vskip 0.2truecm
\par
\noindent
$X_1=X$,
\vskip 0.2truecm
\par
\noindent
$X_2=\left(a_1,~...,~a_{j-1},~
\underbrace{
a_j+\sqrt{1+\sum_{k \in \{1,...,n\}\setminus\{j\}}\limits b_k^2}+b_j \cdot \ii}_{j-{\rm th}~{\rm coordinate}},~
a_{j+1},~...,~a_n\right)$,
\vskip 0.2truecm
\par
\noindent
$X_3=S+\left(0,~...,~0,\underbrace{~~~~~~\ii~~~~~~}_{j-{\rm th}~{\rm coordinate}},0,~...,~0\right)=$
\newpage
\begin{eqnarray*}
\Biggl(a_{1}+\sqrt{\frac{1+(b_{j}-1)^{2}}{n-1}},~...,~
a_{j-1}+\sqrt{\frac{1+(b_{j}-1)^2}{n-1}}~,& &\\
\underbrace{
a_{j}+\sqrt{1+\sum_{{k \in \{1,...,n\}\setminus\{j\}}}\limits b_k^2}~+~\ii}
_{j-{\rm th}~{\rm coordinate}}~,& &\\
a_{j+1}+\sqrt{\frac{1+(b_{j}-1)^2}{n-1}},~...,~
a_{n}+\sqrt{\frac{1+(b_{j}-1)^2}{n-1}}\Biggr),
\end{eqnarray*}
\vskip 0.2truecm
\par
\noindent
$X_4=P_1+(0,~...,~0,\underbrace{~~~~~~\ii~~~~~~}_{j-{\rm th}~{\rm coordinate}},0,~...,~0)$,
\vskip 0.2truecm
\par
\noindent
$X_5=P_2+(0,~...,~0,\underbrace{~~~~~~\ii~~~~~~}_{j-{\rm th}~{\rm coordinate}},0,~...,~0)$,
\vskip 0.2truecm
\par
\noindent
. . . . . . . . . . . . . . . . . . . . . . . . . . .
\vskip 0.2truecm
\par
\noindent
$X_{m+3}=P_m+(0,~...,~0,\underbrace{~~~~~~\ii~~~~~~}_{j-{\rm th}~{\rm coordinate}},0,~...,~0)$,
\vskip 0.2truecm
\par
\noindent
$X_{m+4}=T+(0,~...,~0,\underbrace{~~~~~~\ii~~~~~~}_{j-{\rm th}~{\rm coordinate}},0,~...,~0)
=(0,~...,~0,\underbrace{~~~\sqrt{n}+\ii~~~}_{j-{\rm th}~{\rm coordinate}},0,~...,~0)$,
\vskip 0.2truecm
\par
\noindent
$X_{m+5}=(\underbrace{\ii,...,\ii}_{n-{\rm times}})$
\vskip 0.2truecm
\par
\noindent
belong to ${\C}^n$ and satisfy:
\vskip 0.2truecm
\par
\noindent
$\varphi_n(X_{k-1},X_k)=1$ for each $k \in \{2,3,...,m+5\}$,
$\psi_n(X_1,X_2)=b_j^2 \neq 0$, $\psi_n(X_2,X_3)=b_j \neq 0$,
$\psi_n(X_{k-1},X_k)=1$ for each $k \in \{4,5,...,m+5\}$.
\vskip 0.2truecm
\par
\noindent
By Lemma 2 for each $k \in \{2,3,...,m+5\}$
\vskip 0.2truecm
\par
\centerline{$X_k \in \A$ implies $X_{k-1} \in \A$}
\par
\noindent
and
\par
\noindent
\centerline{$X_k \in \B$ implies $X_{k-1} \in \B$.}
\vskip 0.2truecm
\par
\noindent
Therefore, $(\underbrace{\ii,...,\ii}_{n-{\rm times}})=X_{m+5} \in \A$ implies $X=X_1 \in \A$, and also,
$(\underbrace{\ii,...,\ii}_{n-{\rm times}})=X_{m+5} \in \B$ implies $X=X_1 \in \B$.
\vskip 0.2truecm
\par
\noindent
{\bf Theorem 2.} If $n \geq 2$, $f:{\C}^n \to {\C}^n$
preserves unit distance and $f_{|{\R}^{\scriptstyle n}}={\rm id}({\R}^n)$, then
$f={\rm id}({\C}^n)$ or $f=(\underbrace{\tau,...,\tau}_{n-{\rm times}})$.
\vskip 0.2truecm
\par
\noindent
{\it Proof.} By Lemma 4
\vskip 0.2truecm
\par
\noindent
\centerline{$(\underbrace{\ii,...,\ii}_{n-{\rm times}}) \in \A$~~implies~~${\C}^n \setminus {\R}^n \subseteq \A$}
\par
\noindent
and
\par
\noindent
\centerline{$(\underbrace{\ii,...,\ii}_{n-{\rm times}}) \in \B$~~implies~~${\C}^n \setminus {\R}^n \subseteq \B$.}
\vskip 0.2truecm
\par
\noindent
Obviously, ${\R}^n \subseteq \A$ and ${\R}^n \subseteq \B$. Therefore,
\vskip 0.2truecm
\par
\noindent
\centerline{$\A={\C}^n$ and $f={\rm id}({\C}^n)$,~~if $(\underbrace{\ii,...,\ii}_{n-{\rm times}}) \in \A$,}
\par
\noindent
and also,
\par
\noindent
\centerline{$\B={\C}^n$ and $f=(\underbrace{\tau,...,\tau}_{n-{\rm times}})$,~~if $(\underbrace{\ii,...,\ii}_{n-{\rm times}}) \in \B$.}
\vskip 0.4truecm
\par
Theorem 2 has a simpler proof under the additional assumption that $f$ is continuous.
We need a topological lemma.
\vskip 0.2truecm
\par
\noindent
{\bf Lemma 5.} ${\C}^n \setminus {\R}^n$ is connected for each $n \geq 2$.
\vskip 0.2truecm
\par
\noindent
{\it Proof.} Let $X=(a_1+b_1 \cdot \ii,~...,~a_n+b_n \cdot \ii) \in {\C}^n \setminus {\R}^n$, where
$a_1,...,a_n,b_1,...,b_n \in \R$. We choose $j \in \{1,...,n\}$
with $b_j \neq 0$. Then
$$Y:=\left(0,~...,~0,\underbrace{~~~\frac{b_j}{|b_j|} \cdot \ii~~~}_{j-{\rm th}~{\rm coordinate}}, 0,~...,~0\right) \in {\C}^n \setminus {\R}^n$$
and the segments $(\underbrace{\ii,...,\ii}_{n-{\rm times}})Y$ and $YX$ are disjoint from ${\R}^n$.
These segments form a path joining $(\underbrace{\ii,...,\ii}_{n-{\rm times}})$ 
and $X$, $X$ is an abitrary point in ${\C}^n \setminus {\R}^n$. It proves that ${\C}^n \setminus {\R}^n$
is connected.
\vskip 0.4truecm
\par
Since ${\rm id}({\C}^n \setminus {\R}^n)$ and
$(\underbrace{\tau,...,\tau}_{n-{\rm times}})_{|{\C}^{\scriptstyle n} \setminus {\R}^{\scriptstyle n}}$ are continuous,
for continuous $f$ in Theorem 2
\vskip 0.2truecm
\par
\noindent
~~~~~~~~~~~~~~~~~~~~~~~~~~$\A \setminus {\R}^n=\{X \in {\C}^n \setminus {\R}^n: f(X)=X\}$
\par
\noindent
and
\par
\noindent
~~~~~~~~~~~~~~~~~~~~~~~~~~$\B \setminus {\R}^n=\{X \in {\C}^n \setminus {\R}^n: f(X)=(\underbrace{\tau,...,\tau}_{n-{\rm times}})(X)\}$
\vskip 0.2truecm
\par
\noindent
are the closed subsets of ${\C}^n \setminus {\R}^n$. Obviously,
\vskip 0.2truecm
\par
\noindent
\centerline{$(\A\setminus {\R}^n) \cup (\B \setminus {\R}^n)={\C}^n\setminus {\R}^n$}
\par
\noindent
and
\par
\noindent
\centerline{$(\A \setminus {\R}^n) \cap (\B \setminus {\R}^n)=\emptyset$.}
\par
\noindent
Hence by Lemma 5
\vskip 0.1truecm
\par
\noindent
~~~~~~~~~~~~~~~~~~~~~~~~~~$\A \setminus {\R}^n={\C}^n \setminus {\R}^n$~~~or~~~$\B \setminus {\R}^n={\C}^n \setminus {\R}^n$.
\par
\noindent
Thus,
\par
\noindent
~~~~~~~~~~~~~~~~~~~~~~~~~~$\A={\C}^n~~~{\rm and}~~~f={\rm id}({\C}^n)$
\par
\noindent
or
\par
\noindent
~~~~~~~~~~~~~~~~~~~~~~~~~~$\B={\C}^n~~~{\rm and}~~~f=(\underbrace{\tau,...,\tau}_{n-{\rm times}})$.
\vskip 0.2truecm
\par
\noindent
{\bf Theorem 3.} If $n \geq 2$ and $f:{\C}^n \to {\C}^n$
preserves all positive distances, then $f$ has a form $I \circ (\underbrace{\rho,...,\rho}_{n-{\rm times}})$,
where $\rho \in \{{\rm id}(\C),\tau\}$ and
$I:{\C}^n \to {\C}^n$ is an affine mapping with orthogonal linear part.
\vskip 0.2truecm
\par
\noindent
{\it Proof.} By {\bf (4)} there exists an affine
$I:{\C}^n \to {\C}^n$ such that a linear part of $I$ is orthogonal
and $I_{|{\R}^{\scriptstyle n}}=f_{|{\R}^{\scriptstyle n}}$.
Then $I^{-1} \circ f:{\C}^n \to {\C}^n$ preserves all positive
distances, $(I^{-1} \circ f)_{|{\R}^{\scriptstyle n}}={\rm id}({\R}^n)$.
By Theorem 2~~$I^{-1} \circ f={\rm id}({\C}^n)$ or $I^{-1} \circ f=(\underbrace{\tau,...,\tau}_{n-{\rm times}})$.
In the first case $f=I \circ (\underbrace{{\rm id}(\C),...,{\rm id}(\C)}_{n-{\rm times}})$,
in the second case $f=I \circ (\underbrace{\tau,...,\tau}_{n-{\rm times}})$.
\vskip 0.2truecm
\par
As a corollary of Lemma 1 and Theorem 3 we get:
\vskip 0.2truecm
\par
\noindent
{\bf Theorem 4.} If $n \geq 2$ and a continuous $f:{\C}^n \to {\C}^n$
preserves unit distance, then $f$ has a form $I \circ (\underbrace{\rho,...,\rho}_{n-{\rm times}})$,
where $\rho \in \{{\rm id}(\C),\tau\}$ and $I:{\C}^n \to {\C}^n$
is an affine mapping with orthogonal linear part.
\vskip 0.4truecm
\par
Any bijective $f:{\C}^n \to {\C}^n$~$(n \geq 3)$
that preserves unit distance has a form $I \circ (\rho,...,\rho)$,
where $\rho: \C \to \C$ is a field isomorphism and
$I:{\C}^n \to {\C}^n$ is an affine mapping with orthogonal
linear part; it follows from Theorem 2 in \cite{Schroder}.
\vskip 0.2truecm
\par
The author proved in \cite{Tyszka4}:
\vskip 0.2truecm
\par
\noindent
{\bf (13)}~~~each unit-distance preserving mapping from
${\C}^2$ to ${\C}^2$ has a form $I \circ (\rho,\rho)$, where
$\rho:\C \to \C$ is a field homomorphism and
$I:{\C}^2 \to {\C}^2$ is an affine mapping with orthogonal linear part.
\vskip 0.2truecm
\par
\noindent
The first proof of {\bf (13)} in \cite{Tyszka4} is based on
the results of \cite{Tyszka1} and \cite{Tyszka3}.
The second proof of {\bf (13)} in \cite{Tyszka4} is based on
the result of \cite{Schaeffer}. Obviously, for $n=2$ Theorem~3 follows
from~{\bf (13)}.
\vskip 0.2truecm
\par
If a continuous $\sigma: \C \to \C$ is a field homomorphism,
then $(\sigma,\sigma):{\C}^2 \to {\C}^2$ is continuous and preserves unit distance, also
$(\sigma,\sigma)((0,0))=(0,0)$, $(\sigma,\sigma)((1,0))=(1,0)$, $(\sigma,\sigma)((0,1))=(0,1)$.
Therefore, by Theorem 4 $\sigma={\rm id}(\C)$ or $\sigma=\tau$.
We have obtained an alternative geometric proof of a well-known result:
\vskip 0.2truecm
\par
\noindent
{\bf (14)}~~the only continuous field endomorphisms of
$\C$ are ${\rm id}(\C)$ and $\tau$.
\vskip 0.2truecm
\par
\noindent
An algebraic proof of {\bf (14)} may be found in [6, Lemma 1, p.~356].
Conversely, for $n=2$ Theorem 4 follows from {\bf (13)} and {\bf (14)}.
\newpage
\par
Let $d \in \C \setminus \{0\}$ and
$\C \ni x \stackrel{{\tau}_{\scriptstyle d}}{\longrightarrow} \frac{d}{\tau(d)} \cdot \tau(x) \in \C$.
As a consequence of Theorem~4 we get:
\vskip 0.2truecm
\par
\noindent
{\bf Theorem 5.} If $n \geq 2$ and a continuous $f:{\C}^n \to {\C}^n$
preserves distance $d \in \C \setminus \{0\}$, then $f$ has a form
$I \circ (\underbrace{\rho,...,\rho}_{n-{\rm times}})$,
where $\rho \in \{{\rm id}(\C),\tau_d\}$ and
$I:{\C}^n \to {\C}^n$ is an affine mapping with
orthogonal linear part.
\vskip 0.2truecm
\par
\noindent
{\bf Corollary.} If $n \geq 2$, $d_1,d_2 \in \C \setminus \{0\}$,
$\frac{d_{\scriptstyle 1}^{\scriptstyle 2}}{d_{\scriptstyle 2}^{\scriptstyle 2}} \not\in \R$ and
a continuous $f:{\C}^n \to {\C}^n$ preserves distances $d_1$ and $d_2$,
then $f$ is an affine mapping with orthogonal linear part.
\vskip 0.2truecm
\par
Theorems 1--4 do not hold for $n=1$ because the mapping
$\C \ni z \longrightarrow z+{\rm Im}(z) \in \C$
preserves all real distances. In case $n=1$, there is an easy result:
\vskip 0.2truecm
\par
\noindent
{\bf Theorem 6.} Let $f: \C \to \C$ and for each $x,y \in \C$~~$\varphi_1(x,y) \in \R$ implies
$\varphi_1(x,y)=\varphi_1(f(x),f(y))$. Then $f$ has a form $I \circ \rho$,
where $\rho \in \{{\rm id}(\C),\tau\}$ and
$I: \C \to \C$ is an affine mapping with orthogonal linear part.

Apoloniusz Tyszka\\
Technical Faculty\\
Hugo Ko\l{}\l{}\k{a}taj University\\
Balicka 104, 30-149 Krak\'ow, Poland\\
E-mail address: {\it rttyszka@cyf-kr.edu.pl}
\end{document}